# BETA POISSON-G FAMILY OF DISTRIBUTIONS: ITS PROPERTIES AND APPLICATION WITH FAILURE TIME DATA

**Laba Handique*[a], Subrata Chakraborty [b], Farrukh Jamal [c]**

[a] Department of Statistics Madhabdev University, Assam-784164, India
[b] Department of Statistics Dibrugarh University, Assam-786004, India
[c] Department of Statistics Govt. S.A P/G College, Bahawalpur, Pakistan

**\*Corresponding Author  Email: handiquelaba@gmail.com          *Preprint (May 19, 2020)***

## ABSTRACT

A new generalization of the family of Poisson-G is called beta Poisson-G family of distribution. Useful expansions of the probability density function and the cumulative distribution function of the proposed family are derived and seen as infinite mixtures of the Poisson-G distribution. Moment generating function, power moments, entropy, quantile function, skewness and kurtosis are investigated. Numerical computation of moments, skewness, kurtosis and entropy are tabulated for select parameter values. Furthermore, estimation by methods of maximum likelihood is discussed. A simulation study is carried at under varying sample size to assess the performance of this model. Finally suitability check of the proposed model in comparison to its recently introduced models is carried out by considering two real life data sets modeling.

_________________________

**KEYWORDS**: Beta generated family; Poisson-G family; Maximum likelihood; AIC; A; W

## 1. Introduction

Let $r(t)$ be the probability density function (pdf) of a random variable $T \in (\alpha, \beta)$ for $-\infty \le \alpha < \beta < \infty$ and let $W[G(x)]$ be a function of the cumulative distribution function (cdf) of a random variable $X$ such that $W[G(x)]$ satisfies the following conditions:

$$\begin{cases} (i) \quad W[G(x)] \in [\alpha, \beta] \\ (ii)\, W[G(x)] \text{ is differenti able and monotonica lly non - decreasing , and} \\ (iii)\, W[G(x)] \to \alpha \text{ as } x \to -\infty \text{ and } W[G(x)] \to \beta \text{ as } x \to \infty \end{cases} \quad (1)$$

Alzaatreh *et al*. (2013) defined the T-X family cdf by

$$F(x) = \int_{\alpha}^{W[G(x)]} r(t)\,dt = R\{W[G(x)]\} \quad (2)$$

where $W[G(x)]$ satisfies the conditions (1). The pdf corresponding to (2) is given by

$$f(x) = \left\{ \frac{d}{dx} W[G(x)] \right\} r\{W[G(x)]\}$$

where $G(x)$ is the cdf of any baseline distribution.

The Poisson-G family of distribution with cdf is given by (see Kumaraswamy Poisson-G family, Chakraborty *et al*., 2019c)

$$F^{\mathrm{PG}}(x;\lambda) = \frac{1 - e^{-\lambda G(x)}}{1 - e^{-\lambda}}, \ \lambda \in R - \{0\}; \ n = 1, 2, \ldots \tag{3}$$

The corresponding pdf and hrf Poisson-G family is given by

$$f^{\mathrm{PG}}(x;\lambda) = \frac{\lambda g(x) e^{-\lambda G(x)}}{1 - e^{-\lambda}}, \ \lambda \in R - \{0\}; \ -\infty < x < \infty \tag{4}$$

and

$$h^{\mathrm{PG}}(x) = \frac{\lambda g(x) e^{-\lambda G(x)}}{e^{-\lambda G(x)} - e^{-\lambda}}$$

Of course the characterizations for positive and negative values of the additional parameter will be different.

In this paper, we propose a new class of continuous distributions called the ***beta generalized Poisson-G*** ($BP-G(m,n,\lambda)$) family by taking $W(G(x)) = (1 - e^{-\lambda})^{-1}(1 - e^{-\lambda G(x)})$ and $r(t) = (B(m,n))^{-1} t^{m-1}(1-t)^{n-1}$, $0 < t < 1$ and $m,n,\alpha > 0$. The sf and cdf of $BP-G(m,n,\lambda)$ family of distribution is given by

$$\overline{F}^{BP-G}(x;m,n,\lambda) = 1 - I_{\frac{1-e^{-\lambda G(x)}}{1-e^{-\lambda}}}(m,n), \ F^{BP-G}(x;m,n,\lambda) = I_{\frac{1-e^{-\lambda G(x)}}{1-e^{-\lambda}}}(m,n) \cdot \tag{5}$$

where $B(m,n)$ is the beta function and $I_x(a,b) = B(m,n)^{-1} \int_0^x t^{m-1}(1-t)^{n-1} dt$ the incomplete beta function ratio, $G(x)$ is the baseline cdf. The corresponding pdf and hrf is given by

$$f^{BP-G}(x;m,n,\lambda) = \frac{\lambda g(x) e^{-\lambda G(x)} [1 - e^{-\lambda G(x)}]^{m-1} [e^{-\lambda G(x)} - e^{-\lambda}]^{n-1}}{B(m,n)(1 - e^{-\lambda})^{m+n-1}} \tag{6}$$

and $$h^{BP-G}(x;m,n,\lambda) = \frac{\lambda g(x) e^{-\lambda G(x)} [1 - e^{-\lambda G(x)}]^{m-1} [e^{-\lambda G(x)} - e^{-\lambda}]^{n-1}}{B(m,n)(1 - e^{-\lambda})^{m+n-1}[1 - I_{\frac{1-e^{-\lambda G(x)}}{1-e^{-\lambda}}}(m,n)]}$$

where $m, n, \lambda > 0$, $x > 0$ and $g(x) = G'(x)$ is the baseline distribution.

Some of the notable distributions derived the beta Marshall-Olkin-G family (Alizadeh *et al*., 2015), beta generated Kumaraswamy-G family (Handique *et al*., 2017a), beta Kumaraswamy Marshall-Olkin-G family (Handique and Chakraborty, 2017a), beta generalized Marshall-Olkin Kumaraswamy-G family (Handique and Chakraborty, 2017b), beta Marshall-Olkin Kumaraswamy-G family (Chakraborty *et al*., 2018), Marshall-Olkin Kumaraswamy-G family of distribution (Handique *et al*., 2017b), generalized Marshall-Olkin Kumaraswamy-*G* family of distribution (Chakraborty and Handique, 2017), Kumaraswamy generalized Marshall-Olkin-G family (Chakraborty and Handique, 2018), exponentiated generalized Marshall-Olkin-G family (Handique *et al*., 2018), Odd moment exponential-G family (Haq *et al*., 2018), Zografos-Balakrishnan-Burr-XII family (Emrah *et al*., 2018), generalized Marshall-Olkin Burr-XII family (Handique and Chakraborty, 2018), odd log-logistic Burr X family (Rana *et al*., 2018), Zografos-Balakrishnan-Frechet distribution (Chakraborty *et al*., 2019a), exponentiated generalized extended Gompertz distribution (Thiago *et al*., 2019), zero truncated Poisson distribution (Abouelmagd *et al*., 2019), odd log-logistic Burr III family (Handique *et al.,* 2019a), generalized Marshall-Olkin Burr-III distribution (Chakraborty *et al*., 2019b), generalized modified exponential-G family (Handique *et al*., 2019b) and new T-X family of distribution (Handique *et al*., 2020) among others.

The rest of this article is organized in five more Sections. In Section 2 defined physical basis of new family. In Section 3 we discuss some general results of the family. Different methods of estimation of parameters are presented in Section 4. Application of the proposed family is considered in Section 5. The paper ends with a conclusion in Section 6.

## 2. Physical Basis of BP-G
**Theorem 1:**

If $m$ and $n$ are both integers, then the probability distribution of $BP-G(m,n,\lambda)$ arises as distribution of the $m^{th}$ order statistics of a random sample of size $m+n-1$ from $P-G(\lambda)$ distribution.

**Proof:** Let $T_1, T_2,...,T_{m+n-1}$ be a random sample of size $m+n-1$ from $P-G(\lambda)$ distribution with cdf $\dfrac{1-e^{-\lambda G(x)}}{1-e^{-\lambda}}$. Then the pdf of the $m^{th}$ order statistics $T_{(m)}$ is given by

$$= (m+n-1)!/\{(m-1)!\ [(m+n-1)-m]\ !\}\left\{\frac{1-e^{-\lambda G(x)}}{1-e^{-\lambda}}\right\}^{m-1}\left\{\frac{e^{-\lambda G(x)}-e^{-\lambda}}{1-e^{-\lambda}}\right\}^{(m+n-1)-m}$$

$$\times\frac{\lambda\,g(x)\,e^{-\lambda G(x)}}{1-e^{-\lambda}}$$

$$= \Gamma(m+n)/\{\Gamma(m)\,\Gamma(n)\}\frac{\lambda\,g(x)e^{-\lambda G(x)}[1-e^{-\lambda G(x)}]^{m-1}[e^{-\lambda G(x)}-e^{-\lambda}]^{n-1}}{(1-e^{-\lambda})^{m+n-1}}$$

$$=\frac{\lambda\,g(x)e^{-\lambda G(x)}[1-e^{-\lambda G(x)}]^{m-1}[e^{-\lambda G(x)}-e^{-\lambda}]^{n-1}}{B(m,n)(1-e^{-\lambda})^{m+n-1}}$$

## 3. General results
In this Section we derive some general results for the proposed $BP-G(m,n,\lambda)$ family.

### 3.1 Expansions of pdf and cdf

By using binomial expansion in (6), we obtain

$$f^{BP-G}(x;m,n,\lambda) = f^{PG}(x;\lambda)\sum_{j=0}^{n-1}\mu_j F^{PG}(x;\lambda)^{j+m-1} \tag{7}$$

$$= \sum_{j=0}^{n-1}\mu'_j \frac{d}{dx}[F^{PG}(x;\lambda)]^{j+m} \tag{8}$$

where, $\mu'_j = \dfrac{(-1)^j}{B(m,n)(j+m)}\binom{n-1}{j}$ and $\mu_j = \mu'_j(j+m)$

The density function of $BP-G(m,n,\lambda)$ can be expressed as a mixture of $P\text{-}G(\lambda)$ density and cdf.

Now to expand the cdf $F^{BP-G}(x;m,n,\lambda) = I_{\frac{1-e^{-\lambda G(x)}}{1-e^{-\lambda}}}(m,n)$ we use the result to get

$$I_z(a,b) = \frac{B_z(a,b)}{B(a,b)} = \frac{z^a}{B(a,b)}\sum_{i=0}^{\infty}\binom{b-1}{i}\frac{(-1)^i}{(a+i)}z^i$$

Remembering that $B(z;a,b) = B_z(a,b) = z^a\sum_{i=0}^{\infty}\frac{(1-b)_i}{i!(a+i)}z^i$, where $(x)_i$ is a Pochhammer

symbol.

$$= z^a\sum_{i=0}^{\infty}\frac{(-1)^i(b-1)!}{i!(b-i-1)!(a+i)}z^i = z^a\sum_{i=0}^{\infty}\binom{b-1}{i}\frac{(-1)^i}{(a+i)}z^i$$

(See "Incomplete Beta Function" *From Math World*-A Wolfram Web Resource. http://mathworld. Wolfram.com/Incomplete Beta Function. html).

$$F^{BP-G}(x;m,n,\lambda) = \frac{1}{B(m,n)}\left(\frac{1-e^{-\lambda G(x)}}{1-e^{-\lambda}}\right)^m\sum_{p=0}^{\infty}\binom{n-1}{p}\frac{(-1)^p}{(m+p)}\left(\frac{1-e^{-\lambda G(x)}}{1-e^{-\lambda}}\right)^p$$

$$= \sum_{p=q=0}^{\infty}\sum_{r=0}^{q}\frac{(-1)^{p+q+r}}{B(m,n)(m+p)}\binom{n-1}{p}\binom{m+p}{q}\binom{q}{r}F^{PG}(x;\lambda)^r$$

Now in the summation exchanging the indices $q$ and $r$ in the sum symbol, we get

$$F^{BP-G}(x;m,n,\lambda) = \sum_{p=r=0}^{\infty}\sum_{q=r}^{\infty}\frac{(-1)^{p+q+r}}{B(m,n)(m+p)}\binom{n-1}{p}\binom{m+p}{q}\binom{q}{r}F^{PG}(x;\lambda)^r$$

and then $F^{BP-G}(x;m,n,\lambda) = \sum_{r=0}^{\infty}\psi_r F^{PG}(x;\lambda)^r$ \hfill (9)

where $\psi_r = \sum_{p=0}^{\infty} \sum_{q=r}^{\infty} \dfrac{(-1)^{p+q+r}}{B(m,n)(m+p)} \dbinom{n-1}{p} \dbinom{m+p}{q} \dbinom{q}{r}$.

## 3.2 Probability Weighted Moments

The probability weighted moments (PWM), first proposed by Greenwood *et al*. (1979), are expectations of certain functions of a random variable whose mean exists. The $(p,q,r)^{th}$ PWM of $T$ is defined by $\quad \Gamma_{p,q,r} = \int_{-\infty}^{\infty} x^p F(x)^q \, [1-F(x)]^r \, f(x)\, dx$.

From equation (7) the $s^{th}$ moment of $T$ can be written as

$$E(X^s) = \int_0^{\infty} x^s \, f^{BP-G}(x;m,n,\lambda) \, dx = \sum_{j=0}^{n-1} \mu_j \int_0^{\infty} x^s \, F^{PG}(x;\lambda)^{j+m-1} f^{PG}(x;\lambda) \, dx$$

$$= \sum_{j=0}^{n-1} \mu_j \, \Gamma_{s,\, j+m-1,\, 0}$$

where $\Gamma_{p,q,r} = \int_0^{\infty} x^p \left( \dfrac{1-e^{-\lambda G(x)}}{1-e^{-\lambda}} \right)^q \left[ 1 - \dfrac{1-e^{-\lambda G(x)}}{1-e^{-\lambda}} \right]^r \dfrac{\lambda \, g(x) e^{-\lambda G(x)}}{1-e^{-\lambda}} \, dx$ is the PWM of

$P$-$G(\lambda)$ distribution. Therefore the moments of the $BP-G(m,n,\lambda)$ may be expressed in terms of the PWMs of $P$-$G(\lambda)$.

## 3.3 Moment Generating Function

The moment generating function of $BP-G(m,n,\lambda)$ family can be easily expressed in terms of those of the exponentiated $P$-$G(\lambda)$ distribution using the results of Section 3.1. For example using equation (8) it can be seen that

$$M_X(s) = E[e^{sX}] = \int_0^{\infty} e^{sx} f^{BP-G}(x;m,n,\lambda) \, dx = \int_0^{\infty} e^{sx} \sum_{j=0}^{n-1} \mu_j' \frac{d}{dx}[F^{PG}(x;\lambda)]^{j+m} \, dx$$

$$= \sum_{j=0}^{n-1} \mu_j' \int_0^{\infty} e^{sx} \frac{d}{dx}[F^{PG}(x;\lambda)]^{j+m} \, dx = \sum_{j=0}^{n-1} \mu_j \, M_X(s)$$

where $M_X(s)$ is the mgf of a $P$-$G(\lambda)$ distribution.

Mean, variance, skewness and kurtosis are computed for some values of parameters and tabulated in Table 1.

**Table 1**: Mean, variance, skewness and kurtosis of the $BP-E(m,n,\lambda,\beta)$ distribution with different values of $m,n,\lambda$ and $\beta$

| $m$ | $n$ | $\lambda$ | $\beta$ | Mean | Variance | Skewness | Kurtosis |
|---|---|---|---|---|---|---|---|
| 5 | 1 | 2 | 2 | 0.1712 | 0.1562 | 3.1602 | 16.1856 |
| 4 | 1 | 2 | 2 | 0.2054 | 0.1677 | 2.8544 | 14.1588 |
| 3 | 1 | 2 | 2 | 0.2366 | 0.1683 | 2.6804 | 13.3436 |
| 5 | 2 | 2 | 2 | 0.0715 | 0.0350 | 3.3612 | 17.3340 |
| 5 | 3 | 2 | 2 | 0.0378 | 0.0130 | 3.7480 | 20.2250 |
| 5 | 4 | 2 | 2 | 0.0220 | 0.0059 | 4.2707 | 119.7388 |
| 5 | 4 | 3 | 2 | 0.0724 | 0.0101 | 1.5336 | 0.8385 |
| 5 | 4 | 4 | 2 | 0.0925 | 0.0062 | 0.9379 | 21.4520 |
| 5 | 4 | 5 | 2 | 0.0884 | 0.0032 | 0.9365 | 104.2939 |
| 8 | 4 | 5 | 1 | 0.2354 | 0.0213 | 0.5179 | 4.3436 |
| 8 | 4 | 5 | 2 | 0.4708 | 0.0855 | 0.5179 | 4.3436 |
| 2 | 2 | 1.5 | 2 | 0.0566 | 0.0216 | 3.8989 | 24.2507 |
| 2 | 2 | 3 | 3 | 0.0784 | 0.0087 | 2.4670 | 14.3546 |
| 2 | 3 | 5 | 5 | 0.0241 | 0.0004 | 0.1225 | 46.4268 |
| 2 | 3 | 5 | 6 | 0.0201 | 0.0002 | 2.9764 | 27.5143 |
| 3 | 3 | 5 | 6 | 0.0271 | 0.0004 | -0.7494 | 49.5413 |
| 4 | 3 | 5 | 6 | 0.0330 | 0.0005 | -2.0558 | 60.9605 |
| 4 | 4 | 5 | 6 | 0.0253 | 0.0002 | -1.2278 | 30.0848 |
| 4 | 4 | 7.5 | 7.5 | 0.0162 | 0.0001 | -0.0130 | 33.9609 |

### 3.4 Rényi entropy

The Rényi entropy is defined by $I_R(\delta) = (1-\delta)^{-1} \log\left( \int_{-\infty}^{\infty} f(x)^{\delta} dx \right)$ , where $\delta > 0$ and $\delta \neq 1$. Using binomial expansion in equation (6) we can write

$$f^{\text{BPG}}(x;m,n,\lambda)^{\delta} = f^{PG}(x;\lambda)^{\delta} \sum_{i=0}^{\delta(n-1)} \zeta_i \, F^{PG}(x;\lambda)^{i+\delta(m-1)}$$

Thus

$$I_R(\delta) = (1-\delta)^{-1} \log\left( \int_0^{\infty} f^{\text{PG}}(x;\lambda)^{\delta} \sum_{i=0}^{\delta(n-1)} \zeta_i \, [F^{\text{PG}}(x;\lambda)]^{i+\delta(m-1)} dx \right)$$

$$= (1-\delta)^{-1} \log\left( \sum_{i=0}^{\delta(n-1)} \zeta_i \int_0^{\infty} f^{\text{PG}}(x;\lambda)^{\delta} \, [F^{\text{PG}}(x;\lambda)]^{i+\delta(m-1)} dx \right)$$

where

$$\zeta_i = \left( \frac{1}{B(m,n)} \right)^{\delta} \binom{\delta(n-1)}{i} (-1)^m.$$

Numerical computation of this measure is done for some select values of the parameters and reported in table 2.

**Table 2**: Rényi entropy $BP-E(m,n,\lambda,\beta)$ distribution with different parameters values of $m,n,\lambda,\beta$

| Parameters | | | | $\delta$ | | | | | |
|---|---|---|---|---|---|---|---|---|---|
| $m$ | $n$ | $\lambda$ | $\beta$ | **0.2** | **0.5** | **1.5** | **2** | **3** | **5** |
| **5** | **3** | **3** | **2** | -0.1279 | -1.2731 | 1.1414 | 0.3480 | -0.1008 | -0.3683 |
| 4 | 2 | 3 | 2 | 0.3367 | -0.7125 | 0.8227 | 0.2245 | -0.1310 | -0.3546 |
| 2 | 2 | 3 | 2 | 0.2077 | -0.7978 | -0.1493 | -0.5517 | -0.8136 | -0.9927 |
| 2 | 2 | 2 | 2 | 0.2567 | -1.0798 | 1.8670 | 0.9004 | 0.3561 | 0.0354 |
| 2 | 2 | 1 | 1 | 0.6491 | -1.9981 | 8.5351 | 5.6962 | 4.2193 | 3.4339 |
| 1.5 | 2 | 0.5 | 0.5 | 0.8716 | -3.2126 | 14.9539 | 10.1996 | 7.7623 | 6.4947 |

## 3.5 Distribution of order statistic

Suppose $X_1, X_2,...,X_\vartheta$ is a random sample from any $BP-G(m,n,\lambda)$ distribution. Let $X_{u:\vartheta}$ denote the $u^{th}$ order statistic. The pdf of $X_{u:\vartheta}$ can be expressed as

$$f_{u:\vartheta}(x;m,n,\lambda) = \frac{\vartheta!}{(u-1)!(\vartheta-u)!} \sum_{j=0}^{\vartheta-u} (-1)^j \binom{\vartheta-u}{j} f^{BPG}(x) F^{BPG}(x)^{j+u-1} .$$

Now using the general expansion of the pdf and cdf of the $BP-G(m,n,\lambda)$ distribution the pdf of the $u^{th}$ order statistic for of the $BP-G(m,n,\lambda)$ is given as

$$f_{u:\vartheta}(x;m,n,\lambda) = \frac{\vartheta!}{(u-1)!(\vartheta-u)!} \sum_{j=0}^{\vartheta-u} (-1)^j \binom{\vartheta-u}{j} f^{PG}(x;\lambda) \sum_{l=0}^{n-1} \mu_l F^{PG}(x;\lambda)^{l+m-1} .$$

$$\times [\sum_{k=0}^{\infty} \psi_k F^{PG}(x;\lambda)^k]^{j+u-1}$$

where $\mu_l$ and $\psi_k$ defined in section 3.1

$$[\sum_{k=0}^{\infty} \psi_k F^{PG}(x;\lambda)^k]^{j+u-1} = \sum_{k=0}^{\infty} d_{j+u-1,k} F^{PG}(x;\lambda)^k$$

where $\quad d_{j+u-1,k} = \frac{1}{k\mu_0} \sum_{c=1}^{k} [c(j+u)-k]\mu_c d_{j+u-1,\ k-c} \quad$ (Nadarajah *et al.*, 2015)

Therefore the density function of the $u^{th}$ order statistics of $BP-G(m,n,\lambda)$ distribution can be expressed as

$$f_{u:\vartheta}(x;m,n,\lambda) = \frac{\vartheta!}{(u-1)!(\vartheta-u)!} \sum_{j=0}^{\vartheta-u} (-1)^j \binom{\vartheta-u}{j} f^{PG}(x;\lambda) \sum_{l=0}^{n-1} \mu_l F^{PG}(x;\lambda)^{l+m-1} .$$

$$\times \sum_{k=0}^{\infty} d_{j+u-1,k} F^{PG}(x;\lambda)^{k}$$

$$= f^{PG}(x;\lambda) \sum_{l=0}^{n-1} \sum_{k=0}^{\infty} \gamma_{l,k} F^{PG}(x;\lambda)^{k+l+m-1}$$

$$= f^{PG}(x;\lambda) \sum_{l=0}^{n-1} \sum_{k=0}^{\infty} \gamma_{l,k} \sum_{z=0}^{k+l} \binom{k+l}{z} (-1)^{z} [\overline{F}^{PG}(x;\lambda)]^{z}$$

$$= f^{PG}(x;\lambda) \sum_{z=0}^{k+l} \chi_{z} [\overline{F}^{KwG}(t;a,b,\xi)]^{z}$$

$$= -\sum_{z=0}^{k+l} \frac{\chi_{z}}{z+1} \frac{d}{dx} [\overline{F}^{PG}(x;\lambda)]^{z+1}$$

$$= \sum_{z=0}^{k+l} \chi'_{z} f^{PG}(x;\lambda(z+1))$$

where $\quad \chi'_{z} = \sum_{l=0}^{n-1} \sum_{k=0}^{\infty} \dfrac{(-1)^{z+1} \gamma_{l,k}}{z+1} \binom{k+l}{z}$ , $\chi_{z} = \chi'_{z}(z+1)$ ,

$\gamma_{l,k} = \dfrac{\vartheta!}{(u-1)!(\vartheta-u)!} \sum_{j=0}^{\vartheta-u} (-1)^{j} \binom{\vartheta-u}{j} \mu_{l} \, d_{j+u-1,k}$

### 3.6 Quantile function and related results

The quantile function $X$, let $x = Q(u) = F^{-1}(u)$ , can be obtained by inverting (5). Let $z = Q_{m,n}(u)$ be the beta quantile function. Then,

$$x = Q(u) = Q_{G}\left[ -\frac{1}{\lambda} \log\left[1 - Q_{m,n}(u)(1 - e^{-\lambda})\right] \right]$$

For example, let the $G$ be exponential distribution with parameter $\beta > 0$, having pdf and cdf as $g(x:\beta) = \beta e^{-\beta x}$, $x > 0$ and $G(x:\beta) = 1 - e^{-\beta x}$ respectively. Then the $p^{th}$ quantile is obtained as $-(1/\beta) \log[1-p]$. Therefore, the $p^{th}$ quantile $x_{p}$, of $BP-E$ is given by

$$x_{p} = -\frac{1}{\beta} \log\left[ 1 + \frac{1}{\lambda} \log\left[1 - Q_{m,n}(u)(1 - e^{-\lambda})\right] \right]$$

It is possible to obtain an expansion for $Q_{m,n}(u)$ as $z = Q_{m,n}(u) = \sum_{i=0}^{\infty} e_{i} u^{i/m}$ .

(see "Power series" From MathWorld-A Wolfram Web Resource. http://mathworld. wolfram.com / Power Series.html)

where $e_i = [m B(m,n)]^{1/m} d_i$ and $d_0 = 0, d_1 = 1, d_2 = (n-1)/(n+1)$,

$$d_3 = \{(n-1)(m^2 + 3mn - m + 5n - 4)\} / \{2(m+1)^2 (m+2)\},$$

$$d_4 = (n-1)[m^4 + (6n-1)m^3 + (n+2)(8n-5)m^2 + (33n^2 - 30n + 4)m$$

$$+ n(31n - 47) + 18] / [3(m+1)^3 (m+2)(m+3)]\ldots$$

### 3.7 Plots of the skewness and kurtosis

Here the flexibility of skewness and kurtosis of $\mathrm{BP-G}(m,n,\lambda)$ is checked by plotting Galton skewness (S) that measures the degree of the long tail and Morrs (1988) kurtosis (K) that measures the degree of tail heaviness. Theses are respectively defined by

$$S = \frac{Q(6/8) - 2Q(4/8) + Q(2/8)}{Q(6/8) - Q(2/8)} \text{ and } K = \frac{Q(7/8) - Q(5/8) + Q(3/8) - Q(1/8)}{Q(6/8) - Q(2/8)}$$

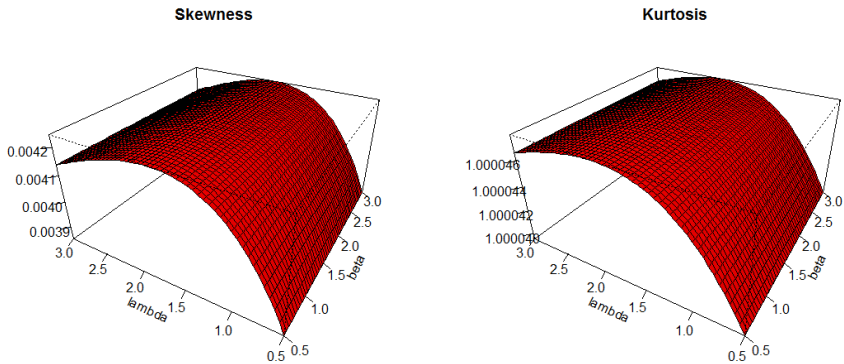

**Fig 1:** Plots of the Galton skewness *S* and Moor kurtosis *K* for the BP-E distribution with parameters
$m = 2,\ n = 3,\ \lambda \in (0.5, 3)$ and $\beta \in (0.5, 3)$

### 3.8 Plots of the pdf and hrf

In this section we have plotted the pdf and hrf of the $BP-G(m,n,\lambda)$ taking $G$ to be Weibull (W) and exponential (E) for some chosen values of the parameters to show the variety of shapes assumed by the family.

The pdf and hrf of these distributions are obtained from $BP-G(m,n,\lambda)$ as follows:

> ➤ The $BP-$Weibull ( $BP-W$ ) distribution

Considering the Weibull distribution (Weibull, 1951) with parameters $\beta > 0$ and $\delta > 0$ having pdf and cdf $g(x) = \delta \beta \, x^{\delta-1} e^{-\beta x^{\delta}}$ and $G(x) = 1 - e^{-\beta x^{\delta}}$ respectively we get the pdf and hrf of $BP\text{-}W\,(m,n,\lambda,\beta,\delta)$ distribution as

$$f^{BP-W}(x) = \frac{\lambda \delta \beta \, x^{\delta-1} e^{-\beta x^{\delta}} \, e^{-\lambda(1-e^{-\beta x^{\delta}})}[1-e^{-\lambda(1-e^{-\beta x^{\delta}})}]^{m-1}[e^{-\lambda(1-e^{-\beta x^{\delta}})}-e^{-\lambda}]^{n-1}}{B(m,n)(1-e^{-\lambda})^{m+n-1}}$$

$$h^{BP-W}(x) = \frac{\lambda \delta \beta \, x^{\delta-1} e^{-\beta x^{\delta}} \, e^{-\lambda(1-e^{-\beta x^{\delta}})}[1-e^{-\lambda(1-e^{-\beta x^{\delta}})}]^{m-1}[e^{-\lambda(1-e^{-\beta x^{\delta}})}-e^{-\lambda}]^{n-1}}{B(m,n)(1-e^{-\lambda})^{m+n-1}[1-I_{\frac{1-e^{-\lambda(1-e^{-\beta x^{\delta}})}}{1-e^{-\lambda}}}(m,n)]}$$

Taking $\delta = 1$ in $BP-W\,(m,n,\lambda,\beta,\delta)$ we get the $BP-E\,(m,n,\lambda,\beta)$ with pdf and hrf is given by

$$f^{BP-E}(x) = \frac{\lambda \beta e^{-\beta x} e^{-\lambda(1-e^{-\beta x})}[1-e^{-\lambda(1-e^{-\beta x})}]^{m-1}[e^{-\lambda(1-e^{-\beta x})}-e^{-\lambda}]^{n-1}}{B(m,n)(1-e^{-\lambda})^{m+n-1}}$$

$$h^{BP-E}(x) = \frac{\lambda \beta e^{-\beta x} e^{-\lambda(1-e^{-\beta x})}[1-e^{-\lambda(1-e^{-\beta x})}]^{m-1}[e^{-\lambda(1-e^{-\beta x})}-e^{-\lambda}]^{n-1}}{B(m,n)(1-e^{-\lambda})^{m+n-1}[1-I_{\frac{1-e^{-\lambda(1-e^{-\beta x})}}{1-e^{-\lambda}}}(m,n)]}$$

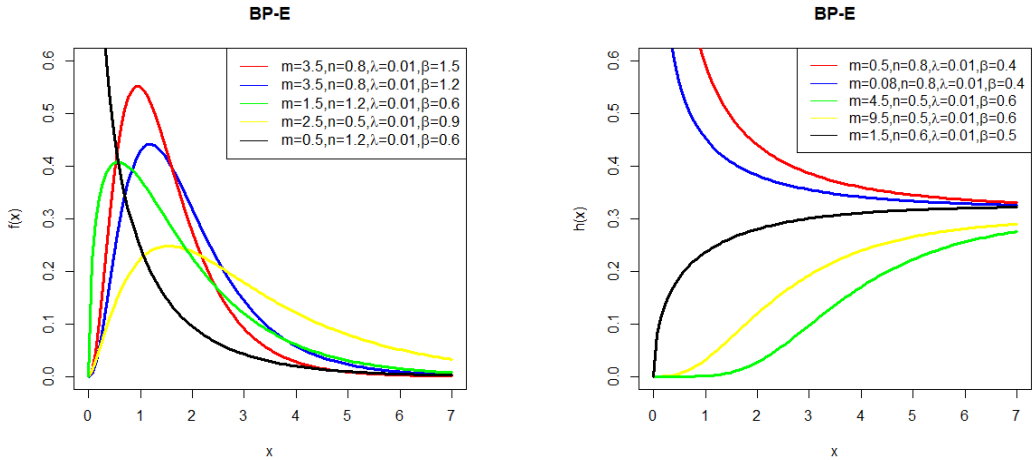

**Fig 2: Density (left) and hazard (right) plots of BP-E**

From the plots in Figure 2 it can be seen that the family is very flexible and can offer many different types of shapes increasing and decreasing failure rate.

## 4. Maximum likelihood estimation

Let $x = (x_1, x_2, ..., x_w)^X$ be a random sample of size $w$ from $BP - E(m, n, \lambda, \beta)$ with parameter vector $\rho = (m, n, \lambda, \beta^T)^T$, where $\beta = (\beta_1, \beta_2, ..., \beta_q)^T$ corresponds to the parameter vector of the baseline distribution $G$. Then the log-likelihood function for $\rho$ is given by

$$\ell = \ell(\rho) = r \log(\lambda) + \sum_{i=1}^{w} \log[g(x_i, \beta)] - \lambda \sum_{i=1}^{w} G(x_i, \beta) + (m-1) \sum_{i=1}^{w} \log(1 - e^{-\lambda G(x_i, \beta)})$$

$$+ (n-1) \sum_{i=1}^{w} \log(e^{-\lambda G(x_i, \beta)} - e^{-\lambda}) - r \log[B(m,n)] - (m+n-1) \sum_{i=1}^{w} \log(1 - e^{-\lambda})$$

The mles are obtained by maximizing the log-likelihood function numerically by using available function from R.

The asymptotic variance-covariance matrix of the MLEs of parameters can obtained by inverting the Fisher information matrix $I(\rho)$ which can be derived using the second partial derivatives of the log-likelihood function with respect to each parameter. The $ij^{th}$ elements of $I_n(\rho)$ are given by

$$I_{ij} = -E[\partial^2 l(\rho)/\partial \rho_i \partial \rho_j], \quad i, j = 1, 2, \cdots, 3+q$$

The exact evaluation of the above expectations may be cumbersome. In practice one can estimate $I_n(\rho)$ by the observed Fisher's information matrix $\hat{I}_n(\hat{\rho}) = (\hat{I}_{ij})$ defined as:

$$\hat{I}_{ij} \approx \left(-\partial^2 l(\rho)/\partial \rho_i \partial \rho_j\right)_{\eta = \hat{\eta}}, \quad i, j = 1, 2, \cdots, 3+q$$

Using the general theory of MLEs under some regularity conditions on the parameters as $n \to \infty$ the asymptotic distribution of $\sqrt{n}(\hat{\rho} - \rho)$ is $N_k(0, V_n)$ where $V_n = (v_{jj}) = I_n^{-1}(\rho)$. The asymptotic behavior remains valid if $V_n$ is replaced by $\hat{V}_n = \hat{I}^{-1}(\hat{\rho})$. Using this result large sample standard errors of $j^{th}$ parameter $\rho_j$ is given by $\sqrt{\hat{v}_{jj}}$.

### 4.1 Simulation

In this section simulation experiments are conducted to evaluate performance of the ML estimates for the $BP - E(m, n, \lambda, \beta)$ distribution with respect to their mean square errors (mse) for sample sizes 50, 100, 200 and 300. The experiment is repeated 3000 times, where in each replication, the ML estimates of the parameters is computed. The average and the MSEs of the estimators is reported in Table 3. From Table 3 its obvious that the estimates are quite stable and more importantly, are close to the true values. More over it is observed that in general the MSEs decreases as $n$ increases.The simulation study therefore shows that the maximum likelihood method is appropriate for estimating the BP-E parameters.

**Table 3**: The average MLEs and MSEs of $BP-E(m,n,\lambda,\beta)$ distribution

| Sample size (w) | True Parameter values | MLE | MSE | True Parameter values | MLE | MSE |
|---|---|---|---|---|---|---|
| 50 | $m=2.2$ | 2.0450 | 0.0598 | $m=2.0$ | 2.0395 | 0.0354 |
| | $n=2.8$ | 2.2570 | 0.1216 | $n=1.8$ | 1.5368 | 0.1119 |
| | $\lambda=0.5$ | 0.4986 | 0.0713 | $\lambda=1.5$ | 1.5068 | 0.0024 |
| | $\beta=2.0$ | 1.0801 | 1.5307 | $\beta=2.0$ | 2.5507 | 1.0447 |
| 100 | $m=2.2$ | 2.0915 | 0.0418 | $m=2.0$ | 2.0160 | 0.0195 |
| | $n=2.8$ | 2.4303 | 0.1078 | $n=1.8$ | 1.6170 | 0.0865 |
| | $\lambda=0.5$ | 0.5088 | 0.0405 | $\lambda=1.5$ | 1.5023 | 0.0018 |
| | $\beta=2.0$ | 1.7711 | 1.0606 | $\beta=2.0$ | 2.1509 | 1.0043 |
| 200 | $m=2.2$ | 2.1903 | 0.0247 | $m=2.0$ | 2.0057 | 0.099 |
| | $n=2.8$ | 2.7899 | 0.0902 | $n=1.8$ | 1.8086 | 0.0511 |
| | $\lambda=0.5$ | 0.5019 | 0.0287 | $\lambda=1.5$ | 1.5010 | 0.0007 |
| | $\beta=2.0$ | 1.9017 | 0.8983 | $\beta=2.0$ | 2.0246 | 0.9005 |
| 300 | $m=2.2$ | 2.2103 | 0.0147 | $m=2.0$ | 1.9998 | 0.0040 |
| | $n=2.8$ | 2.8019 | 0.0512 | $n=1.8$ | 1.8000 | 0.0301 |
| | $\lambda=0.5$ | 0.5002 | 0.0120 | $\lambda=1.5$ | 1.5000 | 0.0004 |
| | $\beta=2.0$ | 2.0017 | 0.4903 | $\beta=2.0$ | 2.0046 | 0.1015 |

Additionally we graphical display of the simulation results in figure 3 and 4. For this we generated for samples of size $w=10$ to $100$ from $BP-E(m,n,\lambda,\beta)$ distribution with true parameters values $m=2.4, n=2.5,\ \lambda=1.5$ and $\beta=1.4,$ and calculated the bias and mean square error (MSE) of the MLEs by repeating the experiment 3000 times.

From the Figures 3 and 4, we observe that when the sample sizes increases, the empirical biases and MSEs approach to zero in all cases justifying the asymptotic normal distribution as an adequate approximation to the finite sample distribution of the MLEs.

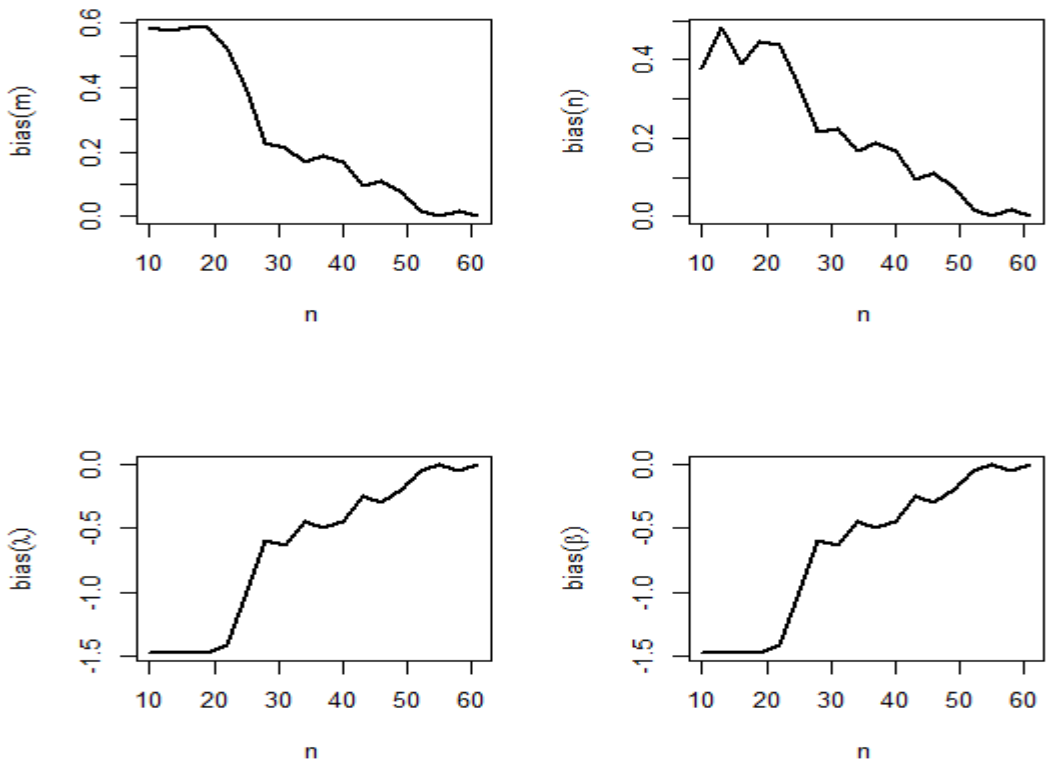

**Fig 3:** The Biases for the parameter values $m = 2.4, n = 2.5, \lambda = 1.5$ and $\beta = 1.4$

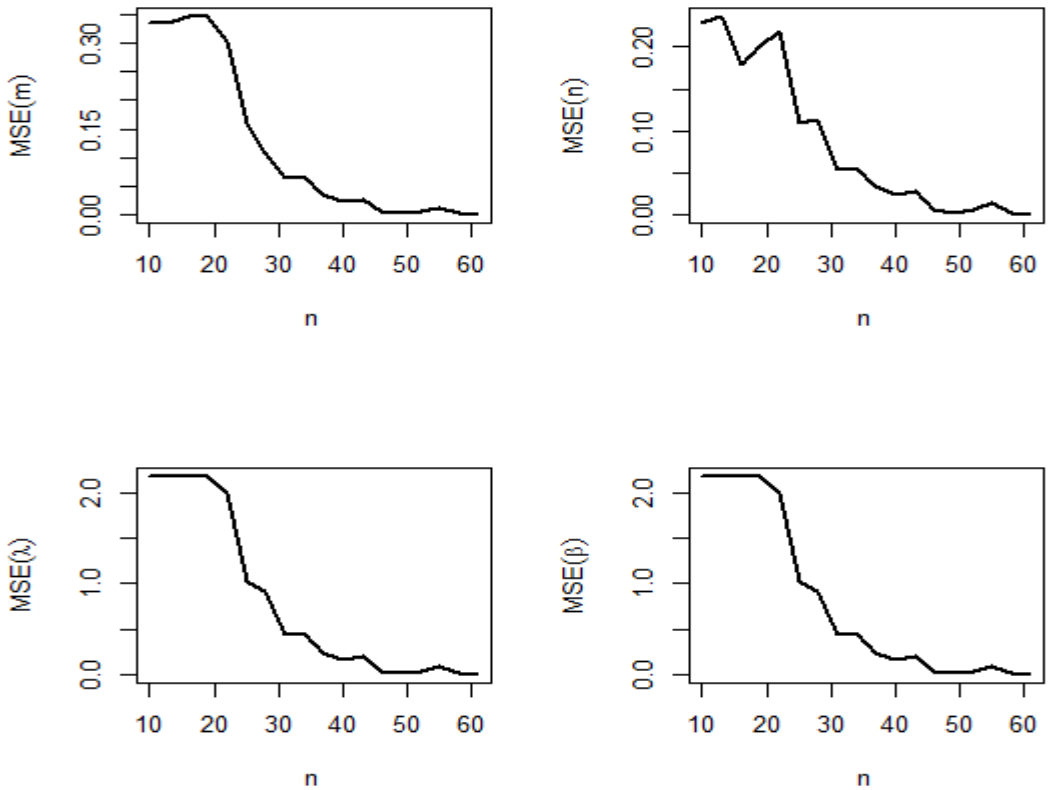

**Fig 4:** The MSEs the parameter values $m = 2.4, n = 2.5, \lambda = 1.5$ and $\beta = 1.4$

## 5. Real life applications

Here we consider fitting of two failure time data sets to show that the distributions from the proposed $BP - E(m, n, \lambda, \beta)$ family can provide better model than the corresponding distributions exponential (Exp), moment exponential (ME), Marshall-Olkin exponential (MO-E), generalized Marshall-Olkin exponential (GMO-E), Kumaraswamy exponential (Kw-E), Beta exponential (BE), Marshall-Olkin Kumaraswamy exponential (MOKw-E) and Kumaraswamy Marshall-Olkin exponential (KwMO-E) distribution.

**Data I:** This data is about survival times (in days) of 72 guinea pigs infected with virulent tubercle bacilli, observed and reported by Bjerkedal (1960), also used by Shibu and Irshad (2016).

**Data II:** The data represents the lifetime data relating to relief times (in minutes) of patients receiving an analgesic. The data was reported by Gross and Clark (1975) and it has twenty (20) observations.

## TTT plots and Descriptive Statistics for the data sets:

The total time on test (TTT) plot (see Aarset, 1987) is a technique to extract the information about the shape of the hazard function. A straight diagonal line indicates constant hazard for the data set, where as a convex (concave) shape implies decreasing (increasing) hazard. The TTT plots for the data sets Fig. 5 indicate that the both data sets have increasing hazard rate.

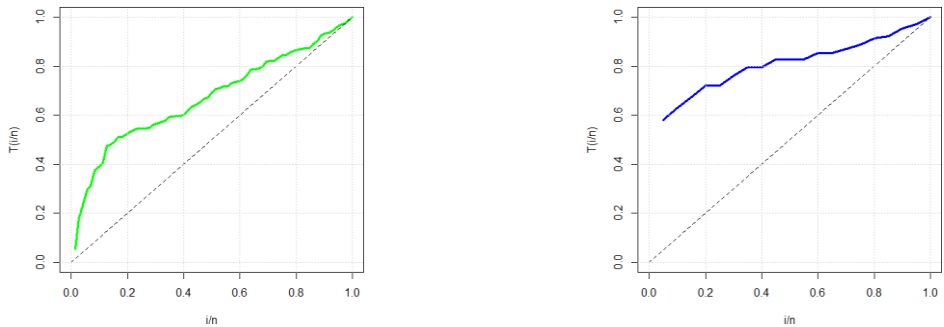

**Fig: 5** TTT-plots for the Data set I and Data set II

We have considered some well known model selection criteria namely the AIC, BIC, CAIC and HQIC and the Kolmogorov-Smirnov (K-S) statistics, Anderson-Darling (A) and Cramer von-mises (W) for goodness of fit to compare the fitted models. We have also provided the asymptotic standard errors and confidence intervals of the mles of the parameters for each competing models. Visual comparison fitted densities and the fitted cdf's are presented in the form of a Histograms and Ogives of the data in Figures 6 and 7. These plots reveal that the proposed distributions provide a good fit to these data.

**Table 4:** Descriptive Statistics for the data sets I and II

| Data Sets | $n$ | Min. | Mean | Median | s.d. | Skewness | Kurtosis | $1^{st}$ Qu. | $3^{rd}$ Qu. | Max. |
|---|---|---|---|---|---|---|---|---|---|---|
| I | 72 | 0.100 | 1.851 | 1.560 | 1.200 | 1.788 | 4.157 | 1.080 | 2.303 | 7.000 |
| II | 20 | 1.100 | 1.900 | 1.700 | 0.704 | 1.592 | 2.346 | 1.475 | 2.050 | 4.100 |

**Table 5 (a):** MLEs, standard errors, confidence intervals (in parentheses) values for the guinea pigs survival time's data set I

| Models | $\hat{\lambda}$ | $\hat{\alpha}$ | $\hat{m}$ | $\hat{n}$ | $\hat{\beta}$ |
|---|---|---|---|---|---|
| Exp $(\beta)$ | --- | --- | --- | --- | 0.540 (0.063) (0.42, 0.66) |
| ME $(\beta)$ | --- | --- | --- | --- | 0.925 (0.077) (0.62, 1.08) |
| MO-E $(\alpha, \beta)$ | --- | 8.778 (3.555) (1.81,15.74) | --- | --- | 1.379 (0.193) (1.00,1.75) |
| GMO-E $(\lambda, \alpha, \beta)$ | 0.179 (0.070) (0.04, 0.32) | 47.635 (44.901) (0, 135.64) | --- | --- | 4.465 (1.327) (1.86. 7.07) |
| Kw-E $(m, n, \beta)$ | --- | --- | 3.304 (1.106) (1.13,5.47) | 1.100 (0.764) (0, 2.59) | 1.037 (0.614) (0, 2.24) |
| B-E $(m, n, \beta)$ | --- | --- | 0.807 (0.696) (0, 2.17) | 3.461 (1.003) (1.49,5.42) | 1.331 (0.855) (0, 3.01) |
| MOKw-E $(\alpha, m, n, \beta)$ | --- | 0.008 (0.002) (0.004, 0.01) | 2.716 (1.316) (0.14, 5.29) | 1.986 (0.784) (0.449, 3.52) | 0.099 (0.048) (0, 0.19) |
| KwMO-E $(\alpha, m, n, \beta)$ | --- | 0.373 (0.136) (0.11, 0.64) | 3.478 (0.861) (1.79, 5.17) | 3.306 (0.779) (1.78, 4.83) | 0.299 (1.112) (0, 2.48) |
| BP-E $(\lambda, m, n, \beta)$ | 0.014 (0.010) (0, 0.03) | --- | 3.595 (1.031) (1.57, 5.62) | 0.724 (1.590) (0, 3.84) | 1.482 (0.516) (0.47, 2.49) |

**Table 5 (b):** Log-likelihood, AIC, BIC, CAIC, HQIC, A, W and KS (*p*-value) values for the guinea pigs survival times data set I

| Models | AIC | BIC | CAIC | HQIC | A | W | KS (*p*-value) |
|---|---|---|---|---|---|---|---|
| Exp $(\beta)$ | 234.63 | 236.91 | 234.68 | 235.54 | 6.53 | 1.25 | 0.27 (0.06) |
| ME $(\beta)$ | 210.40 | 212.68 | 210.45 | 211.30 | 1.52 | 0.25 | 0.14 (0.13) |
| MO-E $(\alpha, \beta)$ | 210.36 | 214.92 | 210.53 | 212.16 | 1.18 | 0.17 | 0.10 (0.43) |
| GMO-E $(\lambda, \alpha, \beta)$ | 210.54 | 217.38 | 210.89 | 213.24 | 1.02 | 0.16 | 0.09 (0.51) |
| Kw-E $(m, n, \beta)$ | 209.42 | 216.24 | 209.77 | 212.12 | 0.74 | 0.11 | 0.08 (0.50) |
| B-E $(m, n, \beta)$ | 207.38 | 214.22 | 207.73 | 210.08 | 0.98 | 0.15 | 0.11 (0.34) |
| MOKw-E $(\alpha, m, n, \beta)$ | 209.44 | 218.56 | 210.04 | 213.04 | 0.79 | 0.12 | 0.10 (0.44) |
| KwMO-E $(\alpha, a, b, \beta)$ | 207.82 | 216.94 | 208.42 | 211.42 | 0.61 | 0.11 | 0.08 (0.73) |
| BP-E $(\lambda, m, n, \beta)$ | 205.42 | 214.50 | 206.02 | 209.02 | 0.55 | 0.08 | 0.09 (0.81) |

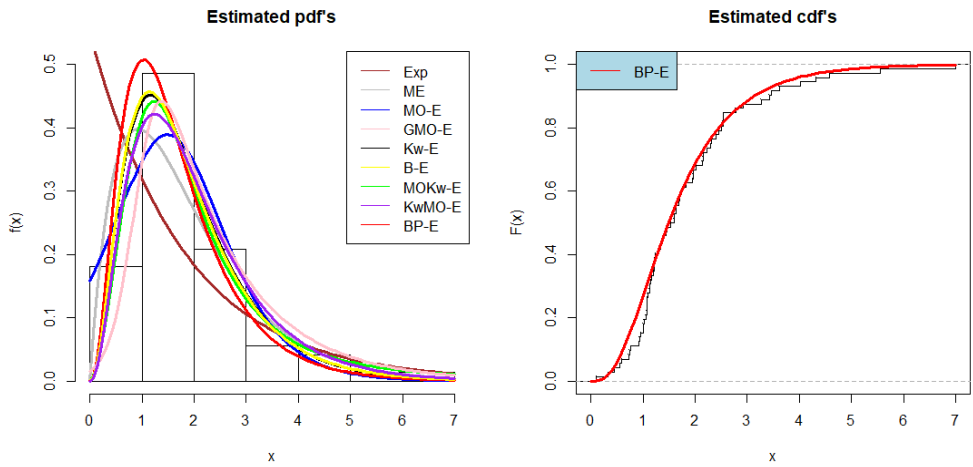

**Fig: 6** Plots of the observed histogram and estimated pdf's on left for the Exp, ME, MO-E, GMO-E, Kw-E, B-E, MOKw-E, KwMO-E, BP-E models and observed ogive and estimated cdf on right for the BP-E data set I.

**Table 6 (a):** MLEs, standard errors, confidence intervals (in parentheses) values for the relief times of patients receiving an analgesic failure time data set II

| Models | $\hat{\lambda}$ | $\hat{\alpha}$ | $\hat{m}$ | $\hat{n}$ | $\hat{\beta}$ |
|---|---|---|---|---|---|
| Exp $(\beta)$ | --- | --- | --- | --- | 0.526 (0.117) (0.29, 0.75) |
| ME $(\beta)$ | --- | --- | --- | --- | 0.950 (0.150) (0.66, 1.24) |
| MO-E $(\alpha, \beta)$ | --- | 54.474 (35.582) (0, 124.21) | --- | --- | 2.316 (0.374) (1.58,3.04) |
| GMO-E $(\lambda, \alpha, \beta)$ | 0.519 (0.256) (0.02, 1.02) | 89.462 (66.278) (0, 219.37) | --- | --- | 3.169 (0.772) (1.66, 4.68) |
| Kw-E $(m, n, \beta)$ | --- | --- | 83.756 (42.361) (0.73, 166.78) | 0.568 (0.326) (0, 1.21) | 3.330 (1.188) (1.00, 5.66) |
| B-E $(m, n, \beta)$ | --- | --- | 81.633 (120.41) (0, 317.63) | 0.542 (0.327) (0, 1.18) | 3.514 (1.410) (0.75, 6.28) |
| MOKw-E $(\alpha, m, n, \beta)$ | --- | 0.133 (0.332) (0, 0.78) | 33.232 (57.837) (0, 146.59) | 0.571 (0.721) (0, 1.98) | 1.669 (1.814) (0, 5.22) |
| KwMO-E $(\alpha, m, n, \beta)$ | --- | 28.868 (9.146) (10.94, 46.79) | 34.826 (22.312) (0, 78.56) | 0.299 (0.239) (0, 0.76) | 4.899 (3.176) (0, 11.12) |
| BP-E $(\lambda, m, n, \beta)$ | 1.965 (0.341) (1.29, 2.63) | --- | 13.396 (1.494) (10.46, 16.32) | 9.600 (1.091) (7.46, 11.73) | 0.244 (0.037) (0.17, 0.32) |

**Table 6 (b):** Log-likelihood, AIC, BIC, CAIC, HQIC, A, W and KS (*p*-value) values for the relief times of patients receiving an analgesic failure time data set II

| Models | AIC | BIC | CAIC | HQIC | A | W | KS (*p*-value) |
|---|---|---|---|---|---|---|---|
| Exp $(\beta)$ | 67.67 | 68.67 | 67.89 | 67.87 | 4.60 | 0.96 | 0.44 (0.004) |
| ME $(\beta)$ | 54.32 | 55.31 | 54.54 | 54.50 | 2.76 | 0.53 | 0.32 (0.07) |
| MO-E $(\alpha, \beta)$ | 43.51 | 45.51 | 44.22 | 43.90 | 0.81 | 0.14 | 0.18 (0.55) |
| GMO-E $(\lambda, \alpha, \beta)$ | 42.75 | 45.74 | 44.25 | 43.34 | 0.51 | 0.08 | 0.15 (0.78) |
| Kw-E $(m, n, \beta)$ | 41.78 | 44.75 | 43.28 | 42.32 | 0.45 | 0.07 | 0.14 (0.86) |
| B-E $(m, n, \beta)$ | 43.48 | 46.45 | 44.98 | 44.02 | 0.70 | 0.12 | 0.16 (0.80) |
| MOKw-E $(\alpha, m, n, \beta)$ | 41.58 | 45.54 | 44.25 | 42.30 | 0.60 | 0.11 | 0.14 (0.87) |
| KwMO-E $(\alpha, m, n, \beta)$ | 42.88 | 46.84 | 45.55 | 43.60 | 1.08 | 0.19 | 0.15 (0.86) |
| BP-E $(\lambda, m, n, \beta)$ | 38.07 | 42.02 | 40.73 | 38.78 | 0.39 | 0.06 | 0.14 (0.91) |

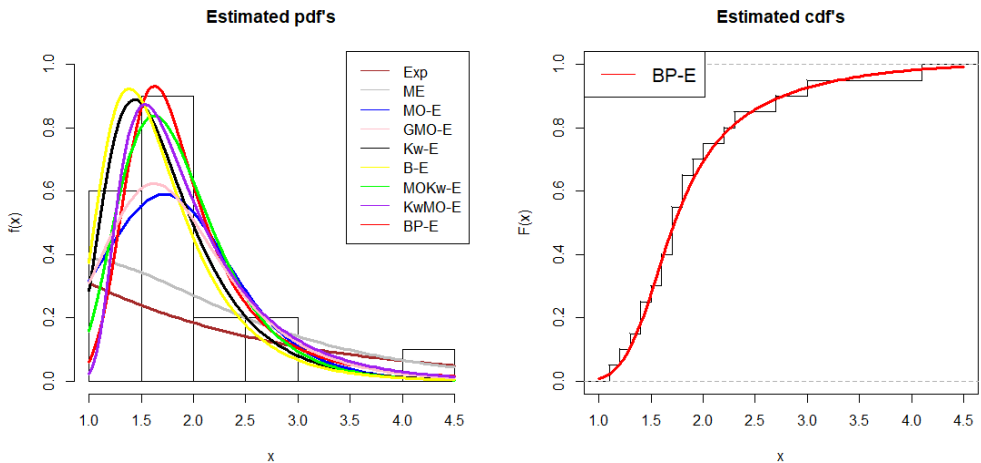

**Fig: 7** Plots of the observed histogram and estimated pdf's on left and observed ogive and estimated cdf's for the Exp, ME, MO-E, GMO-E, Kw-E, B-E, MOKw-E, KwMO-E and BP-E models for data set II.

In the Tables 5(a &b) and 6(a &b) the mle's with standard errors of the parameters for all the fitted models along with their AIC, BIC, CAIC, HQIC, A, W and KS statistic with *p*-value from the fitting results of the data sets I and II are presented respectively.

From the findings presented in the Tables 5(b) and 6(b) on the basis of the lowest value different criteria like AIC, BIC, CAIC, HQIC, A, W and highest *p*-value of the KS statistics the $BP-E$ is found to be a better model than its recently introduced model Exp, ME, MO-E, GMO-E, Kw-E, B-E, MOKw-E and KwMO-E for all the data sets considered here. A visual comparison of the closeness of the fitted densities with the observed histogram and fitted cdf's with the observed ogive of the data sets I and II are presented in the Figures 6 and 7 respectively. These plots also indicate that the proposed distributions provide comparatively closer fit to these data sets.

## 6 Conclusion

A new extension of the Poisson-G family of distributions is introduced and some of its important properties are studied. The MLE for estimating the parameters is discussed. As many as two applications of real life data fitting shows good result in favour of the distributions from the proposed family. Therefore it is expected that the proposed family will be a useful contribution to the existing literature of continuous distributions.